\documentclass[final,3p,times]{elsarticle}
\usepackage{graphics}
\usepackage[utf8]{inputenc} 
\usepackage{latexsym,amssymb,amsthm,amsmath,amsfonts}
\usepackage{url}
\newtheorem{theorem}{Theorem}
\newtheorem{proposition}{Proposition}

\newtheorem{definition}{Definition}
\newtheorem{remark}{Remark}

\newcommand{\R}{{\mathbb R}}
\newcommand{\B}{{\mathbb B}}

\newcommand{\ux}{\underline{x}}

\newcommand{\uy}{\underline{y}}

\newcommand{\pux}{\partial_{\ux}}
\newcommand{\puy}{\partial_{\uy}}

\newcommand{\C}{\mathbb C}

\newcommand{\A}{{\bf\alpha}}

\begin{document}

\begin{frontmatter}

\title{Constructing steering-type solutions for higher order Cauchy-Riemann equations in $\R^{m+1}$}

\author{Daniel Alfonso Santiesteban$^{(1)}$; Dixan Peña Peña$^{(2)}$; Ricardo Abreu Blaya$^{(1)}$}

\address{
$^1$Facultad de Matem\'aticas, Universidad Aut\'onoma de Guerrero, Mexico\\
$^2$Department of Mathematical Analysis, Faculty of Engineering and Architecture, Ghent University, Belgium}
\ead{danielalfonso950105@gmail.com, dixanpena@gmail.com, rabreublaya@yahoo.es}

\begin{abstract}
The multidimensional Cauchy-Riemann operator provides a framework for studying higher order partial differential equations in $\R^{m+1}$, whose solutions include polymonogenic and polyharmonic functions, among others. In this work, we aim to explicitly construct solutions to such systems, generated from families of complex valued functions which are closed under conjugation and under the action of the complex Cauchy-Riemann operator. Moreover, we prove that precisely some of these solutions also satisfy homogeneous linear differential equations involving the so-called hypercomplex derivative. 
\end{abstract}

\end{frontmatter}
\section{Introduction}

Clifford analysis focuses on the study of the null solutions of the Dirac operator $\pux=\sum_{j=1}^me_j\partial_{x_j}$ or the generalized Cauchy-Riemann operator $\partial_{x_0}+\pux$ in an open region of $\R^m$ or $\R^{m+1}$, respectively. These solutions are known as monogenic functions and represent genuine generalizations of holomorphic functions of the complex plane. Here $(e_1,...,e_m)$ is an orthonormal basis in $\R^m$ underlying the construction of the $2^m$-dimensional universal real Clifford algebra $\R_{0,m}$, and $\ux=\sum_{j=1}^mx_je_j$ is a vector variable defined in $\R^m$. For an in-depth study about Clifford analysis, including historical notes, we refer to e.g. \cite{BDS,DSS,Delang,DK,GHS}. 

Clifford algebras have important applications in various fields of research, including geometry, theoretical physics,  cohomology, $K$-theory, quantum mechanics, computer vision, robotics and digital image processing. They are often used to study quadratic forms and orthogonal groups. These algebras are fundamental for describing geometric symmetries in physical spaces and spacetime, especially in relation to rotations and Lorentz transformations. It is worth mentioning that these algebras provide a mathematical structure that unifies the inner and outer products of vectors, allowing for a more general and elegant treatment of geometric concepts. 

The Clifford algebra $\R_{0,m}$ with generators $e_1,e_2,\dots, e_m$, subject to the basic multiplication rules 
\[
e_i^2=-1,\quad e_ie_{j}=-e_{j}e_i,\quad i,j=1,2,\dots, m,\quad i<j,
\] 
is a real linear associative but non-commutative algebra with identity $1$. We shall freely use the well-known properties of Clifford algebras which the reader can find in many sources. Nevertheless, we compile some of those relevant to our purpose for completeness. The algebra $\R_{0,m}$ can be split up subspace of $k$-vectors $(k\in\mathbb N_m\cup\{0\})$
$$\R_{0,m}^{(k)}=\left\{a\in\R_{0,m}:a=\sum_{|A|=k}a_Ae_A,a_A\in\R\right\},$$
where $A$ runs over all the possible ordered sets $A=\{1\leq j_1<...<j_k\leq m\}$, or $A=\emptyset$, and $e_A:=e_{j_1}e_{j_2}...e_{j_k}$, $e_0=e_{\emptyset}=1$. Namely, $\R_{0,m}=\bigoplus_{k=0}^m\R_{0,m}^{(k)}$. 
Any Clifford number $a\in\R_{0,m}$ may thus be written as
$$a=[a]_0+[a]_1+...+[a]_k+...+[a]_m,$$
where $[]_k$ is the linear projection of $a$ onto the subspace of $k$-vectors $\R_{0,m}^{(k)}$. In particular, $\R_{0,m}^{(1)}$, $\R_{0,m}^{(0)}\oplus\R_{0,m}^{(1)}$ and $\R_{0,m}^{(m)}$ are called, respectively, the space of vectors, paravectors and pseudoscalars in $\R_{0,m}$. Note that $\R^{m+1}$ may be naturally identified with $\R_{0,m}^{(0)}\oplus\R_{0,m}^{(1)}$ by associating to any element $(x_0,x_1,...,x_m)\in\R^{m+1}$ the paravector $X=x_0+\ux=x_0+\sum_{j=1}^mx_je_j$. For a $1$-vector $\underline{v}$ and a $k$-vector $F_k$, their product  $\underline{v}F_k$ splits into a $(k-1)$-vector and a $(k+1)$-vector, namely: 
$$\underline{v}F_k=[\underline{v}F_k]_{k-1}+[\underline{v}F_k]_{k+1},$$
where 
$$[\underline{v}F_k]_{k-1}=\frac{1}{2}[\underline{v}F_k-(-1)^kF_k\underline{v}]$$
and
$$[\underline{v}F_k]_{k+1}=\frac{1}{2}[\underline{v}F_k+(-1)^kF_k\underline{v}].$$
The inner and outer products between $\underline{v}$ and $F_k$ are defined by $\underline{v}\cdot F_k:=[\underline{v}F_k]_{k-1}$ and $\underline{v}\wedge F_k:=[\underline{v}F_k]_{k+1}$, respectively.
  
Conjugation in $\R_{0,m}$ is defined as the anti-involution $a\to\overline{a}$ for which $\overline{e_i}=-e_i$. Therefore
$$\overline{a}=\sum_Aa_A\overline{e}_A,\quad\overline{e}_A=(-1)^{\frac{|A|(|A|+1)}{2}}e_A.$$
A norm $\|.\|$ on $\R_{0,m}$ is defined by $\|a\|^2=[a\overline{a}]_{0}=\sum_{A}a_A^2$.

Throughout this paper, $\Omega$ denotes a non-empty open set of $\R^{m+1}$. We will consider functions $f:\Omega\to\R_{0,m}$ to be written as $f=\sum_{A}f_Ae_A$, where $f_A$ are $\R$-valued functions. The notions of continuity, differentiability and integrability of an $\R_{0,m}$-valued function have the usual component-wise meaning. The space of all $s$-times continuous differentiable functions is denoted by $C^s(\Omega)$. The generalized Cauchy-Riemann operator 
\begin{equation}
\partial_X=\partial_{x_0}+\pux=\partial_{x_0}+\sum_{j=1}^me_j\partial_{x_j}
\end{equation}
factorizes the Laplace operator in $\R^{m+1}$ in the sense that $\Delta=\partial_X\overline{\partial}_X=\overline{\partial}_X\partial_X$. 

\begin{definition}
An $\R_{0,m}$-valued function f defined and continuously differentiable in $\Omega$, is said to be left (resp. right) monogenic in $\Omega$ if $\partial_Xf=0$ (resp. $f\partial_X=0$) in $\Omega$. In a similar fashion is defined monogenicity with respect to $\pux$. Furthermore, functions which are both left and right monogenic are called two-sided monogenic, i.e. functions satisfying the overdetermined system $\partial_Xf=0=f\partial_X$.
\end{definition}

The factorization of the Laplace operator by $\partial_X$ implies that monogenic functions are also harmonic. The null solutions of iterated operator $\partial_X^n$ are referred to as polymonogenic functions or $n$-monogenic functions (see \cite{Ryan}). Since $\overline{\partial}_X^n\partial_X^n=\Delta^{n}$, polymonogenic functions are also polyharmonic. In this way, Clifford analysis may be viewed as a refinement of classical harmonic analysis. Despite the fact that Clifford analysis extends the most essential features of complex analysis, monogenic functions do not share all properties of holomorphic functions. For instance, due to the non-commutativity of the Clifford algebras, the product of monogenic functions is not necessarily a monogenic function. Therefore, it is natural to explore specialized techniques for constructing monogenic functions and, more broadly, for obtaining solutions to higher-order Cauchy-Riemann equations.

There are numerous techniques available to generate monogenic functions. For example, in \cite{DSS} the authors obtain a Laurent expansion into special functions for polyaxially monogenic functions and present an introduction to the theory of separately monogenic functions in $\R^m$. In \cite{dixan}, biaxial monogenic functions are constructed from Funk-Hecke formula combined with Fueter theorem. Tao Qian generalized Fueter's result on the induction of  regular quaternionic functions from holomorphic functions of one complex variable to $\R^{m+1}$ and, therefore, it is shown to be consistent with M. Sce's generalization for $m$ as odd integers \cite{tao}. Frank Sommen defined Cauchy-Riemann operators on an $m$-dimensional smooth surface of $\R^{m+1}$ and improved the Cauchy-Kovalevskaya extension formula \cite{sommen}. Here it is important to mention that, in general, every monogenic function $f(x_0,\ux)$ in $\R^{m+1}$, i.e. $(\partial_{x_0}+\pux)f(x_0,\ux)=0$, is determined by its restriction to the hyperplane $x_0=0$. Conversely, any given real analytic function $f(\ux)$ defined in a region of $\R^m$ has a unique monogenic Cauchy-Kovalevskaya extension $f(x_0,\ux)$.   
Together with the so-called Fischer decomposition, the Cauchy-Kovalevskaya extension is one of the essential ingredients for the construction of
orthonormal bases of homogeneous monogenic polynomials, based on group
theoretical methods and arguments. Hilde De Ridder and Frank Sommen discretize techniques for the construction of axially monogenic functions to the setting of discrete Clifford analysis and consider a Vekua-type system for this
construction \cite{ridder}. Moreover, in \cite{erik} the authors present an analogous of the
class of two-sided axial monogenic functions to the case of axial $\kappa$-hypermonogenic functions.

Of particular interest for our propose is the paper \cite{dixan1}, where the authors introduce a new technique leading to so-called steering monogenic functions. These functions can be roughly described as a class of monogenic functions generated from families of complex valued functions which are closed under conjugation and under the action of the classical Cauchy-Riemann operator. In the present paper we generalize their techniques to address a broader family of higher order Cauchy-Riemann systems. 
\section{Steering two-sided monogenic functions}
Consider the biaxial splitting $\R^{m+1}=\R^2\oplus\R^{m-1}$. In this way, for any $X\in\R^{m+1}$ we may rewrite
$$X=z+\uy,$$
where $z=x_0+x_1e_1$ and $\uy=\sum_{j=2}^{m}x_je_j$.
By the above, we can also split the generalized Cauchy-Riemann operator $\partial_X=\partial_{x_0}+\pux$ as
\begin{equation}
\partial_X=2\partial_{\overline{z}}+\puy, 
\end{equation}
where 
\begin{equation}\label{CRo}
\partial_{\overline{z}}=\frac{1}{2}(\partial_{x_0}+e_1\partial_{x_1})
\end{equation}
and 
\begin{equation}\label{Do2}
\puy=\sum_{j=2}^me_j\partial_{x_j}.
\end{equation}
Note that if the identification $i\to e_1$ is made, then the operator \eqref{CRo} is nothing more than the classical Cauchy-Riemann operator in complex analysis. The operator \eqref{Do2}  can be seen as a Dirac operator in $\R^{m-1}$. 

Let $\R_{0,1}$ denote the real Clifford algebra generated by $e_1$, which is isomorphic to $\C$.  Suppose that $\Phi$ is a family of functions $f(z)$ with values in $\R_{0,1}$ which is closed under conjugation and under the action of the operator $\partial_{\overline{z}}$. That is, for any $f\in\Phi$, $\overline{f}\in\Phi$ and $\partial_{\overline{z}}f$ may be expresed as a linear combination of elements in $\Phi$. We shall consider monogenic functions of the form
\begin{equation}\label{f}
\sum_{j}f_j(z)g_j(X),
\end{equation}
where each $f_j$ belongs to $\Phi$ and each $g_j$ is an $\R_{0,m}$-valued function. As the elements of $\Phi$ "steer" the functions $g_j$ in such a way that \eqref{f} is monogenic, these functions \eqref{f} are called steering monogenic functions (see \cite{dixan1}). In what follows, we show some examples. It should be pointed out that some of the following calculations are essentially the same as those developed in \cite{dixan1}, but for the convenience of the reader it was included here. We tacitly assume that any required
differentiations and convergences are well defined.

Consider 
\begin{equation}\label{exp}
\exp(z)A(X)+\exp(\overline{z})B(X)
\end{equation}  
with $\Phi=\{\exp(z),\exp(\overline{z})\}$. By a straightforward computation we obtain
$$\partial_X[\exp(z)A+\exp(\overline{z})B]=\exp(z)(2\partial_{\overline{z}}A+\puy B)+\exp(\overline{z})(\puy A+2\partial_{\overline{z}}B+2B),$$
whence, if $A$ and $B$ satisfy the system
\begin{equation}\label{rl}
\left\{\begin{array}{rl}
2\partial_{\overline{z}}A+\puy B=0,\\
\puy A+2\partial_{\overline{z}}B+2B=0,
\end{array}\right.
\end{equation}
then \eqref{exp} is left monogenic. In particular, if $A$ and $B$ only depend on the variable $\uy$, then the above system also constitutes a necessary condition for the monogenicity of \eqref{exp}, and it takes the form 
\begin{equation}\label{eq8mon}
\left\{\begin{array}{rl}
\puy B=0,\\
B=-\frac{1}{2}\puy A.
\end{array}\right.
\end{equation}
The first equation in \eqref{eq8mon}  tells us that $B(\uy)$ is left monogenic. Moreover, substitution of the second equation of this system into the first one yields $-\puy^2A=\Delta A=0$, i.e. the function $A(\uy)$ is harmonic. We have actually proved that if $H$ is a harmonic function of $\uy$, then 
\begin{equation}\label{H}
\exp(z)H(\uy)-\frac{1}{2}\exp(\overline{z})(\puy H(\uy))
\end{equation}
is left monogenic. This provides a simple way to generate left monogenic functions.  For example, we can assert that the following functions are left monogenic:
$$\exp(z)x_j -\frac{1}{2}\exp(\overline{z})e_j,$$
$$\exp(z)\left(x_j^2-x_k^2\right)-\exp(\overline{z})\left(x_je_j-x_ke_k\right),$$
$$\exp(z)x_jx_k -\frac{1}{2}\exp(\overline{z})\left(x_ke_j+x_je_k\right),$$
$$\exp(z)\frac{1}{\|\uy\|^{m-3}}+\frac{m-3}{2}\exp(\overline{z})\frac{\uy}{\|\uy\|^{m-1}},$$
where $j,k\geq 2$, $j\not=k$.

Now we analyse the right monogenicity of \eqref{exp}. Firstly, we note that
 \begin{align*}
[\exp(z)A+\exp(\overline{z})B]\partial_X&=\exp(z)A+\exp(z)\partial_{x_0}A+\exp(z)\partial_{x_1}Ae_1+\exp(z)e_1Ae_1\\
&\quad+\exp(\overline{z})B+\exp(\overline{z})\partial_{x_0}B+\exp(\overline{z})(\partial_{x_1}B)e_1-\exp(\overline{z})e_1Be_1\\
&\quad+\exp(z)(A\puy)+\exp(\overline{z})(B\puy)\\
&=\exp(z)(A+2A\partial_{\overline z}+e_1Ae_1+A\puy)+\exp(\overline{z})(B+2B\partial_{\overline{z}}-e_1Be_1+B\puy).
\end{align*}
It is easy to see that therefore if $A$ and $B$ satisfy the system
\begin{equation}\label{rexp}
\left\{\begin{array}{rl}
A+2A\partial_{\overline{z}}+e_1Ae_1+A\puy =0,\\
B+2B\partial_{\overline{z}}-e_1Be_1+B\puy=0,
\end{array}\right.
\end{equation}
then \eqref{exp} is right monogenic. Note the crucial difference between systems \eqref{rl} and \eqref{rexp}: unlike the equations of system \eqref{rl}, the equations of \eqref{rexp} are independent of each other. This means that we will end up with right monogenic functions of the form $\exp(z)A$ and $\exp(\overline{z})B$.

In particular, if $A$ and $B$ only depend on the variable $\uy$, then the system \eqref{rexp} is reduced to the following:
\begin{equation}\label{rexpred}
\left\{\begin{array}{rl}
A+e_1Ae_1+A\puy =0,\\
B-e_1Be_1+B\puy=0.
\end{array}\right.
\end{equation}
\begin{remark}
If $A$ and $B$ take values in the $k$-vector subalgebra generated by $\{e_2,\dots,e_m\}$, then \eqref{rexpred} may be rewritten as
\begin{equation*}
\left\{\begin{array}{rl}
A\puy =-\left(1+(-1)^{k+1}\right)A,\\
B\puy=-\left(1+(-1)^{k}\right)B.
\end{array}\right.
\end{equation*}
\end{remark}
Now we will show that the functions $A$ and $B$ are also harmonic, as in the case of the left monogenicity. We have
\begin{align*}
A\puy^2&=-A\puy+e_1A\puy e_1=A+e_1Ae_1-e_1Ae_1-e_1^2Ae_1^2=0,\\
B\puy^2&=-B\puy-e_1B\puy e_1=B-e_1Be_1+e_1Be_1-e_1^2Be_1^2=0.
\end{align*} 
This also implies that
$$M_1=A+e_1Ae_1,\quad M_2=B-e_1Be_1,$$
are right monogenic functions. Note that $A$ and $B$ cannot be explicitly determined by $M_1$ and $M_2$. 

It is also worth noting that if $A$ is a solution to the first equation of system \eqref{rexpred}, then $B=c\puy A$, $c\in\R$, satisfies the second equation of the system. 

There is another way to provide solutions of \eqref{rexpred}. If $M(\uy)$ is a right monogenic function, then $A=M-e_1Me_1$ is a solution to the first equation of the system and $B=M+e_1Me_1$ is a solution to the second equation. In other words, the functions
$$\exp(z)\left(M(\uy)-e_1M(\uy)e_1\right), \quad\exp(\overline{z})\left(M(\uy)+e_1M(\uy)e_1\right)$$
are right monogenic.

Combining the above results, we may conclude:
\begin{theorem}
A function of the form
$$\exp(z)A(\uy)+\exp(\overline{z})B(\uy)$$
is two-sided monogenic if and only if $A$ is a harmonic satisfying the equation $A+e_1Ae_1+A\puy =0$ and $B$ is a left monogenic function of the form $B=-\frac{1}{2}\puy A$.

In particular, if $M(\uy)$ is a right monogenic function, then 
$$\exp(z)(M(\uy)-e_1M(\uy)e_1)-\frac{1}{2}\exp(\overline{z})(\puy M(\uy)+e_1\puy M(\uy)e_1)$$
is two-sided monogenic. 
\end{theorem}

For example, if we take $M(\uy)=\frac{1}{2}\left(x_2+x_3e_2e_3\right)$, then $\exp(z)(x_2+x_3e_2e_3)-\exp(\overline{z})e_2$ is two-sided monogenic. 

Now consider
\begin{equation}\label{trig}
F(X)=(\cos z)A_1(X)+(\sin z) B_1(X)+(\cos\overline{z}) A_2(X)+(\sin\overline{z}) B_2(X),
\end{equation}
with $\Phi=\{\cos z,\sin z,\cos\overline{z},\sin\overline{z}\}$. A direct computation yields
\begin{align*}
\partial_XF&=(\cos z)(2\partial_{\overline{z}} A_1+\puy A_2)\\&\quad+(\sin z)(2\partial_{\overline{z}} B_1+\puy B_2)\\&\quad+(\cos\overline{z})(\puy A_1+2\partial_{\overline{z}} A_2+2B_2)\\&\quad+(\sin\overline{z})(\puy B_1+2\partial_{\overline{z}} B_2-2A_2).
\end{align*}
Consequently, if
\begin{equation}
\left\{\begin{array}{rl}
2\partial_{\overline{z}}A_1+\puy A_2=0,\\
2\partial_{\overline{z}}B_1+\puy B_2=0,\\
\puy A_1+2\partial_{\overline{z}}A_2+2B_2=0,\\
\puy B_1+2\partial_{\overline{z}}B_2-2A_2=0,
\end{array}\right.
\end{equation}
then \eqref{trig} is left monogenic. Similarly to the exponential steering monogenic functions, if $A_j$ and $B_j$ ($j=1,2$) are functions of $\uy$, then \eqref{trig} is left monogenic if and only if $A_1$ and $B_1$ are harmonic functions and
$$A_2(\uy)=\frac{1}{2}\puy B_1(\uy),\quad B_2(\uy)=-\frac{1}{2}\puy A_1(\uy).$$
This gives rise to left monogenic functions of the form
$$(\cos z)H_1(\uy)+(\sin z)H_2(\uy)+\frac{1}{2}(\cos\overline{z})(\puy H_2(\uy))-\frac{1}{2}(\sin\overline{z})(\puy H_1(\uy)),$$
where $H_1$ and $H_2$ are harmonic functions of $\uy$. 

Let us analyse the right monogenicity of \eqref{trig}. Note that
 \begin{align*}
F\partial_X&=\cos z\left(2A_1\partial_{\overline z}+A_1\puy+B_1+e_1B_1e_1\right)\\
&\quad+\sin z\left(2B_1\partial_{\overline z}+B_1\puy-A_1-e_1A_1e_1\right)\\
&\quad+\cos\overline{z}\left(2A_2\partial_{\overline z}+A_2\puy+B_2-e_1B_2e_1\right)\\
&\quad+\sin\overline{z}\left(2B_2\partial_{\overline z}+B_2\puy-A_2+e_1A_2e_1\right).
\end{align*}
Thus, if 
\begin{equation*}\label{rtrig}
\left\{\begin{array}{rl}
2A_1\partial_{\overline z}+A_1\puy+B_1+e_1B_1e_1=0,\\
2B_1\partial_{\overline z}+B_1\puy-A_1-e_1A_1e_1=0,\\
2A_2\partial_{\overline z}+A_2\puy+B_2-e_1B_2e_1=0,\\
2B_2\partial_{\overline z}+B_2\puy-A_2+e_1A_2e_1=0,
\end{array}\right.
\end{equation*}
then \eqref{trig} is right monogenic. In particular, if $A_j$ and $B_j$ ($j=1,2$) are functions of $\uy$, then \eqref{trig} is right monogenic if and only if
\begin{equation}\label{rtrigd}
\left\{\begin{array}{rl}
A_1\puy+B_1+e_1B_1e_1=0,\\
B_1\puy-A_1-e_1A_1e_1=0,\\
A_2\puy+B_2-e_1B_2e_1=0,\\
B_2\puy-A_2+e_1A_2e_1=0.
\end{array}\right.
\end{equation}

One can prove that $A_j$ and $B_j$ are harmonic and therefore $A_j+(-1)^{j+1}e_1A_je_1$, $B_j+(-1)^{j+1}e_1B_je_1$ are right monogenic. Moreover, it can be readily shown that if $A_1$ and $B_1$ are solutions to the first two equations in \eqref{rtrigd}, then $A_2=c\puy B_1$ and $B_2=-c\puy A_1$, $c\in\R$, satisfy the remaining two equations of the same system.

Note that if $M(\uy)$ and $N(\uy)$ are right monogenic functions, then
$$\cos z\left(M(\uy)-e_1M(\uy)e_1\right)+\sin z\left(N(\uy)-e_1N(\uy)e_1\right),$$
$$\cos\overline{z}\left(M(\uy)+e_1M(\uy)e_1\right)+\sin\overline{z}\left(N(\uy)+e_1N(\uy)e_1\right)$$
are also right monogenic functions.

\begin{theorem}
A function of the form
$$(\cos z)A_1(\uy)+(\sin z)B_1(\uy)+(\cos \overline{z})A_2(\uy)+(\sin \overline{z})B_2(\uy)$$
is two-sided monogenic if and only if $A_1$ and $B_1$ are harmonic satisfying the system 
\begin{equation}\label{rtrigd2}
\left\{\begin{array}{rl}
A_1\puy+B_1+e_1B_1e_1=0,\\
B_1\puy-A_1-e_1A_1e_1=0,
\end{array}\right.
\end{equation}
and $A_2=\frac{1}{2}\puy B_1$ and  $B_2=-\frac{1}{2}\puy A_1$.

In particular, if $M(\uy)$ and $N(\uy)$ are  right monogenic functions, then 
$$(\cos z)(M(\uy)-e_1M(\uy)e_1)+(\sin z)(N(\uy)-e_1N(\uy)e_1)+\frac{1}{2}(\cos\overline{z})(\puy N(\uy)+e_1\puy N(\uy)e_1)-\frac{1}{2}(\sin\overline{z})(\puy M(\uy)+e_1\puy M(\uy)e_1)$$
is two-sided monogenic. 
\end{theorem}


If we now consider power steering left monogenic functions of form 
\begin{equation}\label{powerseries}
A_0(X)+\sum_{k=1}^\infty\left(z^kA_k(X)+\overline{z}^kB_k(X)\right),
\end{equation}
with $\Phi=\{z^k,\overline{z}^k:k\in\mathbb{N}_0\}$, then 
\begin{equation}
\partial_X A_0+2B_1=0
\end{equation}
and for $k\geq 1$
\begin{equation}
\left\{\begin{array}{rl}
2\partial_{\overline{z}}A_k+\puy B_k=0,\\
\puy A_k+2\partial_{\overline{z}}B_k+2(k+1)B_{k+1}=0.
\end{array}\right.
\end{equation}
In particular, if the coefficients in the series \eqref{powerseries} are functions of $\uy$ only, we conclude that $A_k$ ($k\geq 0$) are harmonic functions of $\uy$ and
 \begin{equation}\label{relpower}
 B_k=-\frac{1}{2k}\puy A_{k-1},\quad k\geq 1.
 \end{equation}
This yields left monogenic functions of form
\begin{equation}
H_0(\uy)+\sum_{k=1}^\infty\left(z^kH_k(\uy)-\frac{1}{2k}\overline{z}^k(\puy H_{k-1}(\uy))\right),
\end{equation}
the functions $H_k$ ($k\geq 0$) being harmonic of the variable $\uy$. 

The functions given by \eqref{powerseries} are right monogenic if
\begin{equation}
A_0\partial_X+A_1+e_1A_1e_1+B_1-e_1B_1e_1=0
\end{equation}
and for $k\geq 1$
\begin{equation}
\left\{\begin{array}{rl}
2A_k\partial_{\overline{z}}+A_k\puy+(k+1)A_{k+1}+(k+1)e_1A_{k+1}e_1=0,\\
2B_k\partial_{\overline{z}}+B_k\puy+(k+1)B_{k+1}-(k+1)e_1B_{k+1}e_1=0.
\end{array}\right.
\end{equation}
In particular, if $A_k$ and $B_k$ are functions of $\uy$ only, we have also that $A_k$ and $B_k$ are harmonic functions satisfying
$$A_0\puy+A_1+e_1A_1e_1+B_1-e_1B_1e_1=0$$
and for $k\geq 1$
\begin{equation}
\left\{\begin{array}{rl}
A_k\puy+(k+1)A_{k+1}+(k+1)e_1A_{k+1}e_1=0,\\
B_k\puy+(k+1)B_{k+1}-(k+1)e_1B_{k+1}e_1=0.
\end{array}\right.
\end{equation}
\begin{theorem}
A function of the form
$$A_0(\uy)+\sum_{k=1}^\infty\left(z^kA_k(\uy)+\overline{z}^kB_k(\uy)\right),$$
is two-sided monogenic if and only if the functions $A_k$ are harmonic satisfying the system 
\begin{equation}\label{rtrigd2}
\left\{\begin{array}{rl}
A_0\puy+A_1+e_1A_1e_1-\frac{1}{2}\puy A_0+\frac{1}{2}e_1\puy A_0e_1=0,\\
A_k\puy+(k+1)A_{k+1}+(k+1)e_1A_{k+1}e_1=0,\quad k\geq 1,
\end{array}\right.
\end{equation}
and $B_k=-\frac{1}{2k}\puy A_{k-1}$, $k\geq 1$.

In particular, if $M_k(\uy)$, $k\geq 0$, are  right monogenic functions, then 
$$M_0(\uy)+\sum_{k=1}^\infty\left(z^k(M_k(\uy)-e_1M_k(\uy)e_1)-\frac{1}{2k}\overline{z}^k(\puy M_{k-1}(\uy)+e_1\puy M_{k-1}(\uy)e_1)\right)$$
is two-sided monogenic. 
\end{theorem}
\begin{remark}
We would also like to point out that similar results can be obtained if we consider exponential, trigonometric and power steering monogenic functions of the type
\begin{align}
&A(X)\exp(z)+B(X)\exp(\overline{z}),\\
&A_1(X)\cos z+B_1(X)\sin z+A_2(X)\cos\overline{z}+B_2(X)\sin\overline{z},\\
&A_0(X)+\sum_{k=1}^\infty\left(A_k(X)z^k+ B_{k}(X)\overline{z}^k\right).
\end{align}
This same remark also applies to the results we will obtain in the following sections.
\end{remark}

\section{Exponential steering polymonogenic functions}\label{Section3}
As mentioned in the introduction, a function $f\in C^n(\Omega)$ is $n$-monogenic (left or right) if it satisfies the equation
$$\partial_X^n f=0\quad (\text{respectively},\;f\partial_X^n=0)$$
in $\Omega$.  

Now consider
\begin{equation}\label{F}
F(X)=\exp(z)A(\uy)+\exp(\overline{z})B(\uy),
\end{equation}
where $A$ and $B$ are functions of $\uy$ only. The left-handed iteration of $\partial_X$ takes the form:
\begin{align*}
\partial_X F(X)&=\exp(z)(\puy B)+\exp(\overline{z})(\puy A+2B),\\
\partial_X^2F(X)&=\exp(z)(\puy^2 A+2\puy B)+\exp(\overline{z})(2\puy A+\puy^2B+4B),\\
\partial_X^3 F(X)&=\exp(z)(2\puy^2 A+\puy^3 B+4\puy B)+\exp(\overline{z})(\puy^3A+4\puy^2 B+4\puy A+8B),\\
\partial_X^4 F(X)&=\exp(z)(\puy^4A+4\puy^3B+4\puy^2A+8\puy B)+\exp(\overline{z})(4\puy^3 A+\puy^4 B+12\puy^2B+8\puy A+16B),
\\
\partial_X^5F(X)&=\exp(z)(4\puy^4A+\puy^5 B+12\puy^3B+8\puy^2A+16\puy B)+\exp(\overline{z})(\puy^5 A+6\puy^4B+12\puy^3A+32\puy^2B+16\puy A+32B).\\
&\;\;\vdots
\end{align*}
If $\partial_X^2 F(X)=0$ then it is necessary to satisfy the following system
\begin{equation}
\left\{\begin{array}{rl}
\puy^2A+2\puy B=0,\\
2\puy A+\puy^2B+4B=0.
\end{array}\right.
\end{equation}
It is easy to see that
\begin{align*}
\puy^3 A=-2\puy^2 B=4\puy A+8B,
\end{align*}
and, therefore,
\begin{align*}
\puy^4 A&=4\puy^2 A+8\puy B=0,\\
B&=\frac{1}{8}\puy^3A-\frac{1}{2}\puy A.
\end{align*}
Hence, the functions of form
\begin{equation}\label{sbim}
\exp(z)H(\uy)+\exp(\overline{z})\left(\frac{1}{8}\puy^3 H(\uy)-\frac{1}{2}\puy H(\uy)\right),
\end{equation}
$H(\uy)$ being an arbitrary biharmonic function, are left bimonogenic. 

If now $\partial_X^3 F(X)=0$, then we arrive at the system
\begin{equation}
\left\{\begin{array}{rl}
2\puy^2A+\puy^3 B+4\puy B=0,\\
\puy^3 A+4\puy^2B+4\puy A+8B=0.
\end{array}\right.
\end{equation}
We obtain
\begin{align*}
\puy^4 A&=-8\puy B-4\puy^2 A-4\puy^3 B=-8\puy B-4\puy^2 A+8\puy^2 A+16\puy B=8\puy B+4\puy^2 A,\\
\puy^5 A&=8\puy^2 B+4\puy^3 A=-2\puy^3A-8\puy A-16B+4\puy^3 A=2\puy^3 A-8\puy A-16 B,\\
\puy^6 A&=2\puy^4 A-8\puy^2 A-16\puy B=0.
\end{align*}
This yields left $3$-monogenic functions of form
\begin{equation}\label{trimH}
\exp(z)H(\uy)+\exp(\overline{z})\left(-\frac{1}{16}\puy^5 H(\uy)+\frac{1}{8}\puy^3 H(\uy)-\frac{1}{2}\puy H(\uy)\right),
\end{equation}
where $\Delta^3 H(\uy)=0$. Next, we present the cases $n=4$ and $n=5$, allowing the reader to appreciate the structure and elegance of these iterative systems.

For $\partial_X^4 F(X)=0$, it is necessary that
\begin{equation}
\left\{\begin{array}{rl}
\puy^4 A+4\puy^3 B+4\puy^2 A+8\puy B=0,\\
\puy^4 B+4\puy^3 A+12\puy^2 B+8\puy A+16B=0.
\end{array}\right.
\end{equation}
Note that, in this case,
\begin{align*}
\puy^5 A&=-4\puy^4B-4\puy^3A-8\puy^2B=12\puy^3A+40\puy^2B+32\puy A+64B,\\
\puy^6A&=12\puy^4A+40\puy^3B+32\puy^2 A+64\puy B=-8\puy^3 B-16\puy^2 A-32\puy B,\\
\puy^7 A&=-8\puy^4B-16\puy^3A-32\puy^2B\\&=16\puy^3 A+64\puy^2B+64\puy A+128B\\
&=16\puy^3 A+1.6\puy^5 A-19.2\puy^3 A-51.2\puy A-102.4B+64\puy A+128 B\\
&=-3.2\puy^3 A+1.6\puy^5 A+12.8\puy A+25.6B,\\
\puy^8 A&=16\puy^4 A+64\puy^3B+64\puy^2 A+128\puy B=0.
\end{align*}
Thus, the function of form
\begin{equation}
\exp(z)H(\uy)+\exp(\overline{z})\left(\frac{5}{128}\puy^7H(\uy)-\frac{1}{16}\puy^5 H(\uy)+\frac{1}{8}\puy^3 H(\uy)-\frac{1}{2}\puy H(\uy)\right),
\end{equation}
with $\Delta^4 H(\uy)=0$, is left $4$-monogenic. 

Finally, if $\partial_X^5 F(X)=0$ then
\begin{equation}\label{26}
\left\{\begin{array}{rl}
4\puy^4 A+\puy^5 B+12\puy^3 B+8\puy^2 A+16\puy B=0,\\
\puy^5 A+6\puy^4 B+12\puy^3 A+32\puy^2 B+16\puy A+32B=0.
\end{array}\right.
\end{equation}
Using the equations of \eqref{26} we have
\begin{align*}
\puy^6 A&=-6\puy^5 B-12\puy^4 A-32\puy^3 B-16\puy^2 A-32\puy B=12\puy^4 A+40\puy^3 B+32\puy^2 A+64\puy B,\\
\puy^7 A&=12\puy^5 A+40\puy^4 B+32\puy^3 A+64\puy^2 B=-32\puy^4 B-112\puy^3 A-320\puy^2 B-192\puy A-384B,\\
\puy^8 A&=-32\puy^5 B-112\puy^4 A-320\puy^3 B-192\puy^2 A-384\puy B=16\puy^4 A+64\puy^3 B+64\puy^2 A+128\puy B,\\
\puy^9 A&=16\puy^5 A+64\puy^4 B+64\puy^3 A+128\puy^2 B\\
&=-32\puy^4 B-128\puy^3 A-384\puy^2 B-256\puy A-512 B\\
&=\puy^7 A-16\puy^3 A-64\puy^2 B-64\puy A-128B,\\
\puy^{10} A&=\puy^8A-16\puy^4A-64\puy^3B-64\puy^2A-128\puy B=0.
\end{align*}
Since
$$16\puy^5 A-3\puy^7A-144\puy^3 A-448\puy^2 B-320\puy A-640B=0$$
then
\begin{align*}
7\puy^9 A-10\puy^7A-32\puy^3A+16\puy^5A+128\puy A+256B=0,
\end{align*}
equivalently,
\begin{align*}
B=-\frac{7}{256}\puy^9 A+\frac{5}{128}\puy^7 A-\frac{1}{16}\puy^5A+\frac{1}{8}\puy^3 A-\frac{1}{2}\puy A.
\end{align*}
This gives
\begin{equation}
F(X)=\exp(z)H(\uy)+\exp(\overline{z})\left(-\frac{7}{256}\puy^9H(\uy)+\frac{5}{128}\puy^7H(\uy)-\frac{1}{16}\puy^5 H(\uy)+\frac{1}{8}\puy^3 H(\uy)-\frac{1}{2}\puy H(\uy)\right),
\end{equation}
with $\Delta^5 H(\uy)=0$.

Having analyzed these particular cases, we now turn our attention to the general case.

Let the operator be defined by the following matrix
\begin{equation}\label{19}
 T(A,B)=\begin{pmatrix}
 0&\puy\\
 \puy&2
 \end{pmatrix}\begin{pmatrix}
 A\\B
 \end{pmatrix}.
 \end{equation}
Note that the iterations of $T$ correspond directly to the factors accompanying the expressions $\exp(z)$ and $\exp(\overline{z})$ in $\partial_X^nF(X)$. Using the transformation $\puy\to x$, let us now consider the following linear application
\begin{equation}\label{20}
\widehat{T}(A,B)=\begin{pmatrix}
0&x\\
x&2
\end{pmatrix}\begin{pmatrix}
A\\B
\end{pmatrix}.
\end{equation}
Diagonalizing the matrix associated with $\widehat{T}$ yields
\begin{align*}
&\begin{pmatrix}
0&x\\
x&2
\end{pmatrix}=\begin{pmatrix}
-x&-x\\
-1-\sqrt{1+x^2}&-1+\sqrt{1+x^2}
\end{pmatrix}\begin{pmatrix}
1+\sqrt{1+x^2}&0\\
0&1-\sqrt{1+x^2}
\end{pmatrix}\begin{pmatrix}
-x&-x\\
-1-\sqrt{1+x^2}&-1+\sqrt{1+x^2}
\end{pmatrix}^{-1}.
\end{align*}
Therefore, the $n$-iteration of $\widehat{T}$ is given by

\begin{align*}
&\widehat{T}^n(A,B)=\begin{pmatrix}
0&x\\
x&2
\end{pmatrix}^n\begin{pmatrix}
A\\B
\end{pmatrix}\\
&=\begin{pmatrix}
-x&-x\\
-1-\sqrt{1+x^2}&-1+\sqrt{1+x^2}
\end{pmatrix}\begin{pmatrix}
(1+\sqrt{1+x^2})^n&0\\
0&(1-\sqrt{1+x^2})^n
\end{pmatrix}\begin{pmatrix}
-x&-x\\
-1-\sqrt{1+x^2}&-1+\sqrt{1+x^2}
\end{pmatrix}^{-1}\begin{pmatrix}
A\\B
\end{pmatrix} 
\\
&=\begin{pmatrix}
\frac{1}{2}\frac{x^2}{\sqrt{x^2+1}}[(\sqrt{x^2+1}+1)^{n-1}-((1-\sqrt{x^2+1})^{n-1})]&\frac{1}{2}\frac{x}{\sqrt{x^2+1}}[(\sqrt{x^2+1}+1)^n-(1-\sqrt{x^2+1})^n]\\
\frac{1}{2}\frac{x}{\sqrt{x^2+1}}[(\sqrt{x^2+1}+1)^n-(1-\sqrt{x^2+1})^n]&\frac{1}{2}\frac{1}{\sqrt{x^2+1}}[(\sqrt{x^2+1}+1)^{n+1}-(1-\sqrt{x^2+1})^{n+1}]
\end{pmatrix}\begin{pmatrix}
A\\B
\end{pmatrix}\\
&=\begin{pmatrix}
\sum_{k=0}^{\lfloor (n-2)/2 \rfloor}\sum_{j=0}^k\binom{n-1}{2k+1}\binom{k}{j}x^{2j+2}&\sum_{k=0}^{\lfloor (n-1)/2 \rfloor}\sum_{j=0}^k\binom{n}{2k+1}\binom{k}{j}x^{2j+1}
\\
\sum_{k=0}^{\lfloor (n-1)/2 \rfloor}\sum_{j=0}^k\binom{n}{2k+1}\binom{k}{j}x^{2j+1}&\sum_{k=0}^{\lfloor (n)/2 \rfloor}\sum_{j=0}^k\binom{n+1}{2k+1}\binom{k}{j}x^{2j}
\end{pmatrix}\begin{pmatrix}
A\\B
\end{pmatrix},\quad n\geq 2,
\end{align*} 
where use has been made of the Newton binomial theorem.

So we obtain
\begin{equation}\label{21}
T^n(A,B)=\begin{pmatrix}
\sum_{k=0}^{\lfloor (n-2)/2 \rfloor}\sum_{j=0}^k\binom{n-1}{2k+1}\binom{k}{j}\puy^{2j+2}&\sum_{k=0}^{\lfloor (n-1)/2 \rfloor}\sum_{j=0}^k\binom{n}{2k+1}\binom{k}{j}\puy^{2j+1}
\\
\sum_{k=0}^{\lfloor (n-1)/2 \rfloor}\sum_{j=0}^k\binom{n}{2k+1}\binom{k}{j}\puy^{2j+1}&\sum_{k=0}^{\lfloor (n)/2 \rfloor}\sum_{j=0}^k\binom{n+1}{2k+1}\binom{k}{j}\puy^{2j}
\end{pmatrix}\begin{pmatrix}
A\\B
\end{pmatrix},\quad n\geq 2.
\end{equation}
According to \eqref{21} the function \eqref{F} is left $n$-monogenic if and only if
\begin{equation}\label{30}
\left\{\begin{array}{rl}
\sum_{k=0}^{\lfloor (n-2)/2 \rfloor}\sum_{j=0}^k\binom{n-1}{2k+1}\binom{k}{j}\puy^{2j+2}A+\sum_{k=0}^{\lfloor (n-1)/2 \rfloor}\sum_{j=0}^k\binom{n}{2k+1}\binom{k}{j}\puy^{2j+1}B=0,\\
\sum_{k=0}^{\lfloor (n-1)/2 \rfloor}\sum_{j=0}^k\binom{n}{2k+1}\binom{k}{j}\puy^{2j+1}A+\sum_{k=0}^{\lfloor (n)/2 \rfloor}\sum_{j=0}^k\binom{n+1}{2k+1}\binom{k}{j}\puy^{2j}B=0,
\end{array}\right.
\end{equation}
where here $\puy^0 B:=B$.

Following a repetitive calculation procedure involving both equations of \eqref{30}, analogous to the method applied in the previously discussed particular cases, we arrive at the necessary conclusion:
\begin{equation}\label{33}
\Delta^n A(\uy)=0
\end{equation}
and
\begin{equation}\label{34}
B(\uy)=\sum_{k=1}^nc_k\puy^{2k-1}A(\uy),
\end{equation}
where 
\begin{align}\label{ck}
c_1&=-\frac{1}{2},\nonumber\\
c_k&=-\frac{1}{2^k}\sum_{j=1}^{\lfloor k/2 \rfloor}\sum_{i=1}^j\binom{k+1}{2j+1}\binom{j}{i}c_{k-i},\quad k\geq 2.
\end{align}
It can be directly verified that, by choosing $B$ as defined in relation \eqref{34} and invoking the polyharmonicity of $A$, a solution to \eqref{30} is obtained.

This fact can also be achieved by deduction and induction. We know that every $p$-monogenic function is also $q$-monogenic with $q>p$; however, the opposite is not true. If $F(X)$ is $p$-monogenic then $\Delta^p A(\uy)=0$ and $B(\uy)=\sum_{k=1}^pc_k\puy^{2k-1}A(\uy)$, therefore, $\Delta^q A(\uy)=0$ and since $\puy^{2k-1}A=0$ for $k>p$ then $B(\uy)=\sum_{k=1}^qc_k \puy^{2k-1}A(\uy)=\sum_{k=1}^pc_k\puy^{2k-1}A(\uy)$. That is, the solutions obtained for the $p$-monogenicity of $F(X)$ belongs to the set of solutions for its $q$-monogenicity, when $q>p$. Obviously, the opposite embedding is not satisfied in general. We have the following result:
\begin{theorem}\label{T1}
Let $A$ and $B$ be two functions that depend only on the variable $\uy$. The $\R_{0,m}$-valued function $$F(X)=\exp(z)A(\uy)+\exp(\overline{z})B(\uy)$$ is left $n$-monogenic if and only if $$\Delta^n A(\uy)=0$$ and
$$B(\uy)=\sum_{k=1}^nc_k\puy^{2k-1}A(\uy),$$
where $c_k$ are given by the recursive formulas \eqref{ck}.
\end{theorem}
 \begin{remark}
 Using \eqref{ck}, it can be verified that precisely
 \begin{align*}
c_2&=-\frac{1}{2^2}\sum_{k=1}^1\sum_{j=1}^k\binom{3}{2k+1}\binom{k}{j}c_{2-j}=-\frac{1}{4}\binom{3}{3}\binom{1}{1}c_1=\frac{1}{8},\\
c_3&=-\frac{1}{2^3}\sum_{k=1}^{1}\sum_{j=1}^k\binom{4}{2k+1}\binom{k}{j}c_{3-j}=-\frac{1}{8}\binom{4}{3}\binom{1}{1}c_2=-\frac{1}{16},\\
c_4&=-\frac{1}{2^4}\sum_{k=1}^{2}\sum_{j=1}^k\binom{5}{2k+1}\binom{k}{j}c_{4-j}=-\frac{1}{16}\left(\binom{5}{3}\binom{1}{1}c_{3}+\binom{5}{5}\binom{2}{1}c_3+\binom{5}{5}\binom{2}{2}c_2\right)=-\frac{1}{16}\left(-\frac{10}{16}-\frac{1}{8}+\frac{1}{8}\right)=\frac{5}{128},\\
c_5&=-\frac{1}{2^5}\sum_{k=1}^{2}\sum_{j=1}^k\binom{6}{2k+1}\binom{k}{j}c_{5-j}=-\frac{1}{32}\left(\binom{6}{3}\binom{1}{1}c_{4}+\binom{6}{5}\binom{2}{1}c_{4}+\binom{6}{5}\binom{2}{2}c_{3}\right)=-\frac{1}{32}\left(\frac{100}{128}+\frac{60}{128}-\frac{6}{16}\right)=-\frac{7}{256},
\end{align*}
as we have obtained for the particular cases when $n=\overline{2,5}$.
 \end{remark}
\subsection{Right polymonogenic functions}
Note that
\begin{align*}
F(X)\partial_X&=\exp(z)(A+e_1Ae_1+A\puy)+\exp(\overline{z})(B-e_1Be_1+B\puy),\\
F(X)\partial_X^2&=\exp(z)(2A+2e_1Ae_1+2A\puy+A\puy^2)+\exp(\overline{z})(2B-2e_1Be_1+2B\puy+B\puy^2),\\
F(X)\partial_X^3&=\exp(z)(4A+4e_1Ae_1+4A\puy+3A\puy^2+e_1A\puy^2 e_1 +A\puy^3)\\&\quad+\exp(\overline{z})(4B-4e_1Be_1+4B\puy+3B\puy^2-e_1B\puy^2e_1+B\puy^3).
\end{align*}
The function $F(X)$ given by \eqref{F} is right bimonogenic if and only if
\begin{equation}
\left\{\begin{array}{rl}
2A+2e_1Ae_1+2A\puy+A\puy^2=0,\\
2B-2e_1Be_1+2B\puy+B\puy^2=0.
\end{array}\right.
\end{equation}
Note that if $B=\frac{1}{8}\puy^3 A-\frac{1}{2}\puy A$, then it can be verified that
\begin{align*}
2B-2e_1Be_1+2B\puy+B\puy^2&=\frac{1}{4}\puy^3A-\puy A-\frac{1}{4}e_1\puy^3Ae_1+e_1\puy Ae_1+\frac{1}{4}\puy^3A\puy-\puy A\puy+\frac{1}{8}\puy^3A\puy^2-\frac{1}{2}\puy A\puy^2\\
&=\frac{1}{4}\puy^3A-\puy A-\frac{1}{4}e_1\puy^3Ae_1+e_1\puy Ae_1+\frac{1}{4}\puy^3A\puy-\puy A\puy-\frac{1}{2}\puy A\puy^2\\
&=-\left(\puy A-e_1\puy Ae_1+\puy A\puy+\frac{1}{2}\puy A\puy^2\right)+\frac{1}{4}\puy^3 A-\frac{1}{4}e_1\puy^3 Ae_1+\frac{1}{4}\puy^3A\puy\\
&=0.
\end{align*}
Furthermore,
\begin{align*}
A\puy^3&=-2A\puy+2e_1A\puy e_1-2A\puy^2=-2A\puy+2e_1A\puy e_1+4A+4e_1Ae_1+4A\puy
=2A\puy+2e_1A\puy e_1+4A+4e_1Ae_1,\\
\Delta^2 A&=-2\puy^2A-2e_1\puy^2Ae_1-2\puy^2A\puy=2(4A+4e_1Ae_1+2A\puy+2e_1A\puy e_1+2A\puy-2e_1A\puy e_1+2A\puy^2)=0.
\end{align*}
In other words, the function given by \eqref{sbim} with $H(\uy)$ satisfying $2H(\uy)+2e_1H(\uy)e_1+2H(\uy)\puy+H(\uy)\puy^2=0$, is right bimonogenic. However, there may be other solutions where $B$ does not depend on $A$. 

Recall that $\partial_X^2$ is not a real-valued  operator, in contrast to the squared Dirac operator. As a result, its left and right actions yield different outcomes, and this asymmetry must be accounted for in the analysis.

If $F(X)\partial_X^3=0$ then
\begin{equation}\label{rm3}
\left\{\begin{array}{rl}
4A+4e_1Ae_1+4A\puy+3A\puy^2+e_1A\puy^2e_1+A\puy^3=0,\\
4B-4e_1Be_1+4B\puy+3B\puy^2-e_1B\puy^2e_1+B\puy^3=0.
\end{array}\right.
\end{equation}
By straightforward
calculations, we obtain
\begin{align*}
 A\puy^4&=-4A\puy+4e_1A\puy e_1-4A\puy^2-3A\puy^3+e_1A\puy^3e_1\\
 &=-4A\puy+4e_1A\puy e_1-4A\puy^2-3A\puy^3-4A-4e_1Ae_1-4e_1A\puy e_1-3e_1A\puy^2e_1-A\puy^2\\
 &=8A\puy+4A\puy^2+8A+8e_1Ae_1,\\
 A\puy^5&=8A\puy^2+4A\puy^3+8A\puy-8e_1A\puy e_1\\
 &=8A\puy^2-16A-16e_1Ae_1-16A\puy-12A\puy^2-4e_1A\puy^2e_1+8A\puy-8e_1A\puy e_1\\
 &=-4A\puy^2-16A-16e_1Ae_1-8A\puy-8e_1A\puy e_1-4e_1A\puy^2e_1,\\
 A\puy^6&=-4A\puy^3-16A\puy+16e_1A\puy e_1-8A\puy^2+8e_1A\puy^2e_1+4e_1A\puy^3 e_1\\
 &=16A+16e_1Ae_1+16A\puy+12A\puy^2+4e_1A\puy^2e_1-16A-16e_1Ae_1-16e_1A\puy e_1-12e_1A\puy^2e_1-4A\puy^2\\
 &\quad-16A\puy+16e_1A\puy e_1-8A\puy^2+8e_1A\puy^2e_1\\
 &=0.
\end{align*}
Putting now $B=-\frac{1}{16}\puy^5 A+\frac{1}{8}\puy^3A-\frac{1}{2}\puy A$, it follows that

\begin{align*}
&4B-4e_1Be_1+4B\puy+3B\puy^2-e_1B\puy^2e_1+B\puy^3\\
&=-\frac{1}{4}\puy^5 A+\frac{1}{2}\puy^3A-2\puy A+\frac{1}{4}e_1\puy^5 Ae_1-\frac{1}{2}e_1\puy^3Ae_1+2e_1\puy Ae_1-\frac{1}{4}\puy^5 A\puy+\frac{1}{2}\puy^3A\puy-2\puy A\puy\\
&\quad +\frac{3}{8}\puy^3A\puy^2-\frac{3}{2}\puy A\puy^2-\frac{1}{8}e_1\puy^3A\puy^2e_1+\frac{1}{2}e_1\puy A\puy^2e_1+\frac{1}{8}\puy^3A\puy^3-\frac{1}{2}\puy A\puy^3\\
&=\frac{1}{8}\puy^5 A-\puy^3 A-2\puy A+\frac{1}{8}e_1\puy^5 Ae_1+2e_1\puy Ae_1-\frac{1}{8}\puy^5A\puy-2\puy A\puy\\
&=\puy A\puy+\frac{1}{2}\puy^3 A+\puy A-e_1\puy Ae_1+e_1\puy A\puy e_1+\frac{1}{2}e_1\puy^3 Ae_1+e_1\puy Ae_1-\puy A-\puy^3 A-2\puy A\\
&\quad+2e_1\puy A e_1-\puy^3 A-\frac{1}{2}\puy^3A\puy -\puy A\puy-e_1\puy A\puy e_1-2\puy A\puy\\
&=-2\puy A\puy-2\puy A-\frac{3}{2}\puy^3 A+2e_1\puy Ae_1+\frac{1}{2}e_1\puy^3 A e_1-\frac{1}{2}\puy^3 A\puy\\
&=0.
\end{align*}
Therefore, the function given by \eqref{trimH} with $H(\uy)$ satisfying the first equation in \eqref{rm3}, is right $3$-monogenic. Now we will focus on explaining why this is so. Let
$$U(A,B)=(A+e_1Ae_1+A\puy,B-e_1Be_1+B\puy).$$
The systems we have obtained for the polymonogenicity of $F(X)$ result from the equation
\begin{equation}
U^n(A,B)=0.
\end{equation}
Note that if $A+e_1Ae_1+A\puy=0$ then 
\begin{equation}
A\puy=-A-e_1Ae_1,
\end{equation}
hence,
\begin{equation}\label{rel40}
A\puy^2=-A\puy+e_1A\puy e_1=0.
\end{equation}
It is easy to see that if a function $f$ is of the form $g+e_1ge_1$, where $g$ is another arbitrary function, then $f-e_1fe_1=0$. This argument quickly leads to \eqref{rel40}. Let now $f=A+e_1Ae_1+A\puy$. If $f+e_1fe_1+f\puy=0$ then $f\puy^2=0$ for the previous argument, and thus,
\begin{align}
A\puy^2+e_1A\puy^2e_1+A\puy^2\puy=0,
\end{align}
which implies $A\puy^4=0$. Let now $g=f+e_1fe_1+f\puy$. If $g+e_1ge_1+g\puy=0$ then $g\puy^2=0$, hence $f\puy^4=0$ and finally $A\puy^6=0$. By repeating this procedure, we can affirm that if $U^n(A,B)=0$, then necessarily, $\Delta^n A=0$. The same reasoning can be
applied to $B$ to obtain the same result. On the other hand, note that the first and second components of the iterations of the operator $U$ only differ in the sign of the elements of the form $e_1(.)e_1$. This can be seen if the operator $H(f)=f+f\puy$ is iterated and rewritten in the form $H(f)=\frac{1}{2}(f+e_1fe_1+f\puy+f-e_1fe_1+f\puy)$. Then it is easy to
see that if $B$ is taken as $B=\sum_{j=1}^rc_j\puy^{2j-1}A$, where $c_j$ are arbitrary real constants, then the second of the equations is solved when $U^n(A,B)=0$. Applying $\sum_{j=1}^rc_j\puy^{2j-1}$ on the left side of the first equation of system $U^n(A,B)=0$ yields the second equation of this system, but now instead of $B$ there is $\sum_{j=1}^rc_j\puy^{2j-1}A$.  We therefore obtain the following proposition:
\begin{proposition}
The function 
$$F(X)=\exp(z)A(\uy)+\exp(\overline{z})B(\uy),$$
where
$$\widetilde{U}^n(A)=0,$$
$$B(\uy)=\sum_{k=1}^nc_k\puy^{2k-1}A(\uy),$$
being $\widetilde{U}(A)=A+e_1Ae_1+A\puy$, and $c_k$ are arbitrary real constants, is right $n$-monogenic. In particular, if the constants $c_k$ are taken as in \eqref{ck} then $F(X)$ is also left $n$-monogenic.   
\end{proposition}
\begin{remark}
Although constructing left polymonogenic functions of form \eqref{F} is more complicated than when considering form 
$$A(\uy)\exp(z)+B(\uy)\exp(\overline{z}),$$
we wanted to show here the differences and similarities that arise in these iterative systems.
\end{remark}

\section{Trigonometric steering polymonogenic functions}
Let us now consider
\begin{equation}\label{trigF}
F(X)=(\cos z) A_1(\uy)+(\sin z) B_1(\uy)+(\cos\overline{z}) A_2(\uy)+(\sin\overline{z}) B_2(\uy),
\end{equation}
where $A_j$, $B_j$, $j=1,2$, are functions of $\uy$ only. We have
\begin{align*}
\partial_X F(X)&=(\cos z)(\puy A_2)+(\sin z)(\puy B_2)+(\cos\overline{z})(\puy A_1+2B_2)+(\sin\overline{z})(\puy B_1-2A_2),\\
\partial_X^2 F(X)&=(\cos z)(\puy^2 A_1+2\puy B_2)+(\sin z)(\puy^2 B_1-2\puy A_2)+(\cos\overline{z})(\puy^2 A_2+2\puy B_1-4A_2)+\sin\overline{z}(\puy^2 B_2-2\puy A_1-4B_2),\\
\partial_X^3 F(X)&=(\cos z)(\puy^3 A_2+2\puy^2 B_1-4\puy A_2)+(\sin z)(\puy^3 B_2-2\puy^2 A_1-4\puy B_2)+(\cos\overline{z})(\puy^3 A_1+4\puy^2 B_2-4\puy A_1-8B_2)\\&\quad+(\sin\overline{z})(\puy^3 B_1-4\puy^2 A_2-4\puy B_1+8A_2),\\
\partial_X^4 F(X)&=(\cos z)(\puy^4 A_1+4\puy^3 B_2-4\puy^2 A_1-8\puy B_2)+(\sin z)(\puy^4 B_1-4\puy^3 A_2-4\puy^2 B_1+8\puy A_2)\\&\quad+(\cos\overline{z})(\puy^4 A_2+4\puy^3 B_1-12\puy^2 A_2-8\puy B_1+16A_2)+(\sin\overline{z})(\puy^4 B_2-4\puy^3 A_1-12\puy^2B_2+8\puy A_1+16B_2).
\end{align*}
The function $F(X)$ given by \eqref{trigF} is left bimonogenic if the following system is satisfied:
\begin{equation}\label{rel43}
\left\{\begin{array}{rl}
\puy^2 A_1+2\puy B_2=0,\\
 \puy^2 B_1-2\puy A_2=0,\\
\puy^2 A_2+2\puy B_1-4A_2=0,\\
\puy^2 B_2-2\puy A_1-4B_2=0.
\end{array}\right.
\end{equation} 
Using the equations in \eqref{rel43} we obtain
 \begin{align*}
\puy^3 A_1&=-2\puy^2 B_2=-4\puy A_1-8B_2,\\
\puy^3 B_1&=2\puy^2 A_2=-4\puy B_1+8A_2,\\
\puy^4 A_1&=0,\\
\puy^4 B_1&=0,
\end{align*}
and, therefore,  $A_2=\frac{1}{8}\puy^3 B_1+\frac{1}{2}\puy B_1$ and $B_2=-\frac{1}{8}\puy^3 A_1-\frac{1}{2}\puy A_1$. Whereas function $F(X)$ would be left $3$-monogenic if the following system is satisfied:
\begin{equation}
\left\{\begin{array}{rl}
\puy^3 A_2+2\puy^2 B_1-4\puy A_2=0,\\
 \puy^3 B_2-2\puy^2 A_1-4\puy B_2=0,\\
\puy^3 A_1+4\puy^2 B_2-4\puy A_1-8B_2=0,\\
\puy^3 B_1-4\puy^2 A_2-4\puy B_1+8A_2=0.
\end{array}\right.
\end{equation}
We arrived at the following
\begin{align*}
\puy^4 A_1&=-4\puy^3 B_2+4\puy^2 A_1+8\puy B_2=-4\puy^2 A_1-8\puy B_2,\\
\puy^5 A_1&=-4\puy^3 A_1-8\puy^2 B_2=-2\puy^3 A_1-8\puy A_1-16B_2,\\
\puy^6 A_1&=6\puy^4 A_1+24\puy^2 A_1+48\puy B_2=0,\\
\puy^4 B_1&=4\puy^3 A_2+4\puy^2 B_1-8\puy A_2=-4\puy^2B_1+8\puy A_2,\\
\puy^5 B_1&= -4\puy^3 B_1+8\puy^2 A_2=-2\puy^3 B_1-8\puy B_1+16A_2,\\
\puy^6 B_1&=-2\puy^4 B_1-8\puy^2 B_1+16\puy A_2=0.
\end{align*}
Thus, $A_2=\frac{1}{16}\puy^5 B_1+\frac{1}{8}\puy^3B_1+\frac{1}{2}\puy B_1$ and $B_2=-\frac{1}{16}\puy^5 A_1-\frac{1}{8}\puy^3 A_1-\frac{1}{2}\puy A_1$. Finally, if $\partial_X^4F(X)=0$ then
\begin{equation}
\left\{\begin{array}{rl}
\puy^4 A_1+4\puy^3B_2-4\puy^2A_1-8\puy B_2=0,\\
\puy^4 B_1-4\puy^3 A_2-4\puy^2B_1+8\puy A_2=0,\\
\puy^4 A_2+4\puy^3 B_1-12\puy^2A_2-8\puy B_1+16A_2=0,\\
\puy^4B_2-4\puy^3A_1-12\puy^2B_2+8\puy A_1+16B_2=0.
\end{array}\right.
\end{equation}
Hence
\begin{align*}
\puy^5 A_1&=-4\puy^4 B_2+4\puy^3A_1+8\puy^2B_2=-12\puy^3A_1-40\puy^2B_2+32\puy A_1+64B_2,\\
\puy^6 A_1&=-12\puy^4A_1-40\puy^3B_2+32\puy^2A_1+64\puy B_2=8\puy^3B_2-16\puy^2A_1-32\puy B_2,\\
\puy^7 A_1&=8\puy^4B_2-16\puy^3A_1-32\puy^2B_2=16\puy^3A_1+64\puy^2B_2-64\puy A_1-128B_2\\&=-1.6\puy^5A_1-3.2\puy^3A_1-12.8\puy A_1-25.6B_2,\\
\puy^8A_1&=16\puy^4A_1+64\puy^3B_2-64\puy^2A_1-128\puy B_2=0,\\
\puy^5B_1&=4\puy^4A_2+4\puy^3B_1-8\puy^2A_2=-12\puy^3B_1+40\puy^2A_2+32\puy B_1-64A_2,\\
\puy^6B_1&=-12\puy^4B_1+40\puy^3A_2+32\puy^2B_1-64\puy A_2=-8\puy^3A_2-16\puy^2B_1+32\puy A_2,\\
\puy^7 B_1&=-8\puy^4A_2-16\puy^3B_1+32\puy^2A_2=16\puy^3 B_1-64\puy^2 A_2-64\puy B_1+128A_2\\&=-1.6\puy^5 B_1-3.2\puy^3B_1-12.8\puy B_1+25.6A_2,\\
\puy^8B_1&=16\puy^4B_1-64\puy^3A_2-64\puy^2B_1+128\puy A_2=0.
\end{align*}
Therefore, $A_2=\frac{5}{128}\puy^7B_1+\frac{1}{16}\puy^5B_1+\frac{1}{8}\puy^3B_1+\frac{1}{2}\puy B_1$ and $B_2=-\frac{5}{128}\puy^7A_1-\frac{1}{16}\puy^5A_1-\frac{1}{8}\puy^3A_1-\frac{1}{2}\puy A_1$.

Note that these systems can be interpreted through iterations of the following operator
\begin{equation}
Q(A_1,B_1,A_2,B_2)=\begin{pmatrix}
0&0&\puy&0\\
0&0&0&\puy\\
\puy&0&0&2\\
0&\puy&-2&0
\end{pmatrix}\begin{pmatrix}
A_1\\B_1\\A_2\\B_2
\end{pmatrix}.
\end{equation}
If the transformation $\puy\to x$ is performed, then we focus the analysis on the following linear algebraic application
$$\widehat{Q}(A_1,B_1,A_2,B_2)=\begin{pmatrix}
0&0&x&0\\
0&0&0&x\\
x&0&0&2\\
0&x&-2&0
\end{pmatrix}\begin{pmatrix}
A_1\\B_1\\A_2\\B_2\\
\end{pmatrix}.$$
While diagonalization of the associated matrix remains a viable method, we choose to proceed more efficiently by utilizing the results previously derived in Section \ref{Section3}. So what we do is look at the submatrices
$$\begin{pmatrix}
0&x\\
x&2
\end{pmatrix},\quad\begin{pmatrix}
0&x\\
x&-2
\end{pmatrix}.$$
The $n$th iteration of $\widehat{Q}$, where $n$ is even, will result in a matrix of the form
$$\begin{pmatrix}
a_1&0&0&a_2\\
0&b_1&b_2&0\\
0&-b_2&b_3&0\\
-a_2&0&0&a_3
\end{pmatrix},$$
where
$$\begin{pmatrix}
a_1&a_2\\
-a_2&a_3
\end{pmatrix}=\left(\begin{pmatrix}
0&x\\
x&-2
\end{pmatrix}\begin{pmatrix}
0&x\\
x&2
\end{pmatrix}\right)^{n/2}=\begin{pmatrix}
x^2&2x\\
-2x&x^2-4
\end{pmatrix}^{n/2}$$
and
$$\begin{pmatrix}
b_1&b_2\\
-b_2&b_3
\end{pmatrix}=\left(\begin{pmatrix}
0&x\\
x&2
\end{pmatrix}\begin{pmatrix}
0&x\\
x&-2
\end{pmatrix}\right)^{n/2}=\begin{pmatrix}
x^2&-2x\\
2x&x^2-4
\end{pmatrix}^{n/2}.$$
Whereas if $n$ is odd, then the iterations will result in a matrix of the form
$$\begin{pmatrix}
0&a_1&a_2&0\\
b_1&0&0&b_2\\
b_2&0&0&b_3\\
0&a_2&a_3&0
\end{pmatrix},$$
where
$$\begin{pmatrix}
a_1&a_2\\
a_2&a_3
\end{pmatrix}=\left(\begin{pmatrix}
0&x\\
x&-2
\end{pmatrix}\begin{pmatrix}
0&x\\
x&2
\end{pmatrix}\right)^{(n-1)/2}\begin{pmatrix}
0&x\\x&-2
\end{pmatrix}=\begin{pmatrix}
x^2&2x\\
-2x&x^2-4
\end{pmatrix}^{(n-1)/2}\begin{pmatrix}
0&x\\x&-2
\end{pmatrix}$$
and
$$\begin{pmatrix}
b_1&b_2\\
b_2&b_3
\end{pmatrix}=\left(\begin{pmatrix}
0&x\\
x&2
\end{pmatrix}\begin{pmatrix}
0&x\\
x&-2
\end{pmatrix}\right)^{(n-1)/2}\begin{pmatrix}
0&x\\x&2
\end{pmatrix}=\begin{pmatrix}
x^2&-2x\\
2x&x^2-4
\end{pmatrix}^{(n-1)/2}\begin{pmatrix}
0&x\\x&2
\end{pmatrix}.$$
So, for our purposes, it suffices to diagonalize the matrices
$$\begin{pmatrix}
x^2&2x\\-2x&x^2-4
\end{pmatrix},\quad\begin{pmatrix}
x^2&-2x\\
2x&x^2-4
\end{pmatrix}.$$
We obtain
\begin{align*}
&\begin{pmatrix}
x^2&2x\\-2x&x^2-4
\end{pmatrix}^n\\&=\begin{pmatrix}
1&1\\
-\frac{1}{x}(\sqrt{1-x^2}+1)&\frac{1}{x}(\sqrt{1-x^2}-1)
\end{pmatrix}\begin{pmatrix}
(x^2-2\sqrt{1-x^2}-2)^n&0\\
0&(x^2+2\sqrt{1-x^2}-2)^n
\end{pmatrix}\\&\quad\times\begin{pmatrix}
1&1\\
-\frac{1}{x}(\sqrt{1-x^2}+1)&\frac{1}{x}(\sqrt{1-x^2}-1)
\end{pmatrix}^{-1},
\end{align*}
\begin{align*}
&\begin{pmatrix}
x^2&-2x\\2x&x^2-4
\end{pmatrix}^n\\&=\begin{pmatrix}
1&1\\
\frac{1}{x}(\sqrt{1-x^2}+1)&-\frac{1}{x}(\sqrt{1-x^2}-1)
\end{pmatrix}\begin{pmatrix}
(x^2-2\sqrt{1-x^2}-2)^n&0\\
0&(x^2+2\sqrt{1-x^2}-2)^n
\end{pmatrix}\\&\quad\times\begin{pmatrix}
1&1\\
\frac{1}{x}(\sqrt{1-x^2}+1)&-\frac{1}{x}(\sqrt{1-x^2}-1)
\end{pmatrix}^{-1}.
\end{align*}
Performing a similar analysis and calculations analogous to those in the previous section, we obtain the following theorem:
\begin{theorem}
Let $A_j,B_j$, $j=1,2$, be functions that depend only on the variable $\uy$. The $\R_{0,m}$-valued function $$F(X)=(\cos z) A_1(\uy)+(\sin z) B_1(\uy)+(\cos\overline{z}) A_2(\uy)+(\sin\overline{z}) B_2(\uy)$$ is left $n$-monogenic if and only if $$\Delta^n A_1(\uy)=0,\quad\Delta^n B_1(\uy)=0,$$ 
$$A_2(\uy)=\sum_{k=1}^n(-1)^{k}c_k\puy^{2k-1}B_1(\uy),$$
and
$$B_2(\uy)=\sum_{k=1}^n(-1)^{k+1}c_k\puy^{2k-1}A_1(\uy),$$
where $c_k$ are given by the recursive formulas \eqref{ck}.
\end{theorem}
\section{Power steering polymonogenic functions}
Now consider the series
$$F(X)=A_0(\uy)+\sum_{k=1}^\infty(z^kA_k(\uy)+\overline{z}^kB_k(\uy)),$$
where $A_0$ and $A_k, B_k$, $k\geq 1$, are functions of $\uy$ only. A simple calculation shows that
\begin{align*}
&\partial_X^2 F(X)\\&=\puy^2 A_0+2\puy B_1+2\puy A_1+8B_2+\sum_{k=1}^\infty( z^k(\puy^2 A_k+2(k+1)\puy B_{k+1})+\overline{z}^k(\puy^2B_k+2(k+1)\puy A_{k+1}+4(k+2)(k+1)B_{k+2})).
\end{align*}
If $\partial_X^2 F(X)=0$ then
\begin{equation}\label{rel48}
\puy^2 A_0+2\puy B_1+2\puy A_1+8B_2=0
\end{equation}
and, for all $k\geq 1$,
\begin{equation}\label{rel49}
\left\{\begin{array}{rl}
\puy^2 A_k+2(k+1)\puy B_{k+1}=0,\\
\puy^2 B_k+2(k+1)\puy A_{k+1}+4(k+1)(k+2)B_{k+2}=0.
\end{array}\right.
\end{equation}
Note that
\begin{align*}
\puy^3 A_0&=-2\puy^2 B_1-2\puy^2 A_1-8\puy B_2=-2\puy^2 B_1=8\puy A_2+48B_3,\\
\puy^4 A_0&=0,\\
\puy^3 A_1&=-4\puy^2 B_2=24\puy A_3+192B_4,\\
\puy^4 A_1&=0,\\
\puy^3A_k&=-2(k+1)\puy^2B_{k+1}=2(k+1)[2(k+2)\puy A_{k+2}+4(k+2)(k+3)B_{k+3}],\\
\puy^4 A_k&=4(k+1)(k+2)[\puy^2 A_{k+2}+2(k+3)\puy B_{k+3}]=0.
\end{align*}
Since any solution to the system obtained for when $\partial_XF(X)=0$ must also be a solution for when $\partial_X^2F(X)=0$, then let $B_1=-\frac{1}{2}\puy A_0$ (see relation \eqref{relpower}). By \eqref{rel48} we obtain $B_2=-\frac{1}{4}\puy A_1$. From \eqref{rel49} follows
$$B_k=-\frac{1}{2k}\puy A_{k-1}+\frac{1}{8(k-1)k(k-2)}\puy^3 A_{k-3},\quad k\geq 3.$$
Then, the function 
$$H_0+zH_1-\overline{z}\frac{1}{2}\puy H_0+z^2H_2-\frac{1}{4}\overline{z}^2\puy H_1+\sum_{k=3}^\infty \left[z^kH_k+\overline{z}^k\left(-\frac{1}{2k}\puy H_{k-1}+\frac{1}{8k(k-1)(k-2)}\puy^3H_{k-3}\right)\right],$$
where $H_j$ $(j\geq 0)$ are biharmonic functions in the variable $\uy$, is left bimonogenic.

If now $\partial_X^3F(X)=0$ then necessarily
\begin{equation}\label{rel50}
\puy^3 A_0+4\puy^2 B_1+2\puy^2 A_1+8\puy B_2+8\puy A_2+48B_3=0
\end{equation}
and, for all $k\geq 1$,
\begin{equation}\label{rel51}
\left\{\begin{array}{rl}
\puy^3B_k+2(k+1)\puy^2A_{k+1}+4(k+1)(k+2)\puy B_{k+2}=0,\\
\puy^3A_k+4(k+1)\puy^2B_{k+1}+2(k+1)[2(k+2)\puy A_{k+2}+4(k+2)(k+3)B_{k+3}]=0.
\end{array}\right.
\end{equation}
Based on the argument mentioned above, we can take $B_1=-\frac{1}{2}\puy A_0$ and $B_2=-\frac{1}{4}\puy A_1$ in equation \eqref{rel50} and, thus, $B_3=\frac{1}{48}\puy^3 A_0-\frac{1}{6}\puy A_2$. By the second equation in \eqref{rel51} we obtain that $B_4=-\frac{1}{8}\puy A_3+\frac{1}{192}\puy^3 A_1$. Now, let us calculate the following:
\begin{align*}
\puy^4 A_0&=-4\puy^3 B_1-2\puy^3A_1-8\puy^2B_2-8\puy^2A_2-48\puy B_3=8\puy^2A_2+48\puy B_3-2\puy^3 A_1-8\puy^2 B_2=8\puy^2A_2+48\puy B_3,\\
\puy^5 A_0&=8\puy^3A_2+48\puy^2B_3=-48\puy^2 B_3-384\puy A_4-3840B_5=-\puy^5 A_0+8\puy^3 A_2-384\puy A_4-3840B_5, 
 \\
 \puy^6 A_0&=0,\\
 \puy^4 A_k&=-4(k+1)\puy^3B_{k+1}-4(k+1)(k+2)\puy^2A_{k+2}-8(k+1)(k+2)(k+3)\puy B_{k+3}\\
 &=8(k+1)(k+2)\puy^2 A_{k+2}+16(k+1)(k+2)(k+3)\puy B_{k+3}-4(k+1)(k+2)\puy^2A_{k+2}-8(k+1)(k+2)(k+3)\puy B_{k+3}\\
 &=4(k+1)(k+2)\puy^2 A_{k+2}+8(k+1)(k+2)(k+3)\puy B_{k+3},\\
 \puy^5 A_k&=4(k+1)(k+2)\puy^3 A_{k+2}+8(k+1)(k+2)(k+3)\puy^2 B_{k+3}\\
 &=4(k+1)(k+2)[-2(k+3)\puy^2 B_{k+3}-4(k+3)(k+4)\puy A_{k+4}-8(k+3)(k+4)(k+5)B_{k+5}]
 \\
 &=4(k+1)(k+2)\puy^3 A_{k+2}+2(k+1)(k+2)[-\puy^3 A_{k+2}-4(k+3)(k+4)\puy A_{k+4}-8(k+3)(k+4)(k+5)B_{k+5}],\\
 \puy^6 A_k&=0.
\end{align*}
Therefore, for $k\geq 0$,
$$B_{k+5}=-\frac{1}{16(k+1)(k+2)(k+3)(k+4)(k+5)}\puy^5 A_k+\frac{1}{8(k+3)(k+4)(k+5)}\puy^3 A_{k+2}-\frac{1}{2(k+5)}\puy A_{k+4},$$
or equivalently, for $k\geq 5$
$$B_k=-\frac{1}{16k(k-1)(k-2)(k-3)(k-4)}\puy^5A_{k-5}+\frac{1}{8k(k-1)(k-2)}\puy^3A_{k-3}-\frac{1}{2k}\puy A_{k-1}.$$ 
Following this methodology and in accordance with the reasoning presented in section \ref{Section3}, we obtain the following theorem:
\begin{theorem}
Let $A_0$, $A_j$, $B_j$, $j\geq 1$, be functions that depend only on the variable $\uy$. The $\R_{0,m}$-valued function
$$F(X)=A_0(\uy)+\sum_{k=1}^\infty (z^kA_k(\uy)+\overline{z}^kB_k(\uy)),$$
where 
$$\Delta^n A_k(\uy)=0,\quad k\geq 0,$$
and for $k\geq 1$
$$B_k(\uy)=\sum_{j=1}^{\min(n,\lfloor (k+1)/2\rfloor)} \frac{c_j}{(2j-1)!\binom{k}{k-2j+1}}\puy^{2j-1}A_{k-2j+1}(\uy),$$
with $c_j$ given by the recursive formulas \eqref{ck}, is left $n$-monogenic.
\end{theorem}
\section{Homogeneous linear differential equations involving the hypercomplex derivative}
Recall that the hypercomplex derivative of a monogenic function $f$ is defined as
\begin{equation}\label{HCD}
Df:=\frac{1}{2}(\partial_{x_0}-\pux)f,
\end{equation}
where $\pux=\sum_{j=1}^me_j\partial_{x_j}$ is the Dirac operator in $\R^m$ \cite{gurlhyp}. As a monogenic function $f$ clearly satisfies $\partial_{x_0}f=-\pux f$, it easily follows that 
\begin{equation}
Df=\partial_{x_0}f=-\pux f.
\end{equation}

The exponential steering monogenic functions
\begin{equation}\label{rel54}
F(X)=\exp(z)A(\uy)+\exp(\overline{z})B(\uy),
\end{equation}
where $ A(\uy)$ is harmonic and $B(\uy)=-\frac{1}{2}\puy A(\uy)$, also have an interesting property. First, note that the hypercomplex derivative \eqref{HCD} of a monogenic function can be further rewritten as
\begin{equation}
Df=\partial_{\overline{z}}f+\partial_zf,
\end{equation}
where $\partial_z=\overline{\partial_{\overline{z}}}=\frac{1}{2}(\partial_{x_0}-e_1\partial_{x_1})$. Then the functions of form \eqref{rel54} satisfy
\begin{equation}\label{imDff}
DF=F.
\end{equation} 
This fact reveals that \eqref{rel54} is a nice generalization of the exponential function $\exp(z)$ for the monogenic functions. It is easy to see that if we take $A(\uy)=2x_je_j$, $j\geq 2$, then the left monogenic function
\begin{equation*}
F(X)=2\exp(z)x_je_j+\exp(\overline{z})
\end{equation*}
satisfies \eqref{imDff} and also
\begin{equation*}
F(0)=1,
\end{equation*}
which suggests a greater similarity to the traditional exponential function. Other interesting examples can be constructed.

Now let us use this fact to solve homogeneous linear differential equations with constant coefficients for monogenic functions. That is, to solve equations of the form
\begin{equation}\label{eqdh}
a_0D^nf+a_1D^{n-1}f+...+a_{n-1}Df+a_nf=0,
\end{equation}
where $a_j\in\R$, $j=\overline{0,n}$, and the function $f$ is left monogenic.

Now consider
\begin{equation}\label{Fr}
F_r(X)=\exp(rz)A(\uy)+\exp(r\overline{z})B(\uy),\quad r\in\R.
\end{equation}
Note that
\begin{align*}
\partial_X F_r(X)&=r\exp(rz)A(\uy)+re_1^2\exp(rz)A(\uy)+\exp(r\overline{z})\puy A\\
&\quad+r\exp(r\overline{z})B(\uy)-re_1^2\exp(r\overline{z})B(\uy)+\exp(rz)\puy B\\
&=\exp(rz)\puy B+\exp(r\overline{z})[\puy A+2rB].
\end{align*}
Therefore, $F_r$ is left monogenic if and only if
\begin{align*}
\puy B&=0,\\
\puy A+2rB&=0,
\end{align*}
or equivalently, $\Delta A=0$ and $B=-\frac{1}{2r}\puy A$ when $r\not=0$. Then the exponential steering monogenic functions of the form \eqref{Fr} satisfy
\begin{equation}\label{rel59}
DF_r=rF_r.
\end{equation}
Consider the characteristic equation of \eqref{eqdh}
\begin{equation}\label{eccara}
a_0\lambda^n+a_1\lambda^{n-1}+...+a_{n-1}\lambda +a_n=0.
\end{equation}
Suppose that  $\lambda_1,\lambda_2,...,\lambda_n$ are the roots of equation \eqref{eccara}, among which there may be multiples. If all the roots are real, distinct, and nonzero, then a solution to \eqref{eqdh} can be given by
\begin{equation}
f(X)=\sum_{j=1}^nC_j\left[\exp(\lambda_iz)H_j(\uy)-\frac{1}{2\lambda_j}\exp(\lambda_j\overline{z})\puy H_j(\uy)\right],
\end{equation}
where $H_j$ are harmonic functions and $C_j$ are arbitrary real constants; whereas if one of these roots is zero, say $\lambda_n$, then a solution to \eqref{eqdh} is
\begin{equation}
f(X)=M(\uy)+\sum_{j=1}^{n-1}C_j\left[\exp(\lambda_jz)H_j(\uy)-\frac{1}{2\lambda_j}\exp(\lambda_j\overline{z})\puy H_j(\uy)\right],
\end{equation}
where $M(\uy)$ is an arbitrary left monogenic function.
\begin{proposition}
Suppose that $L_n=a_0D^n+a_1D^{n-1}+...+a_{n-1}D+a_n$ is a partial differential operator with $a_j\in\R$, $j=\overline{0,n}$. If $\lambda_k$ is a real non-zero root of the algebraic equation
$$a_0\lambda^n+a_1\lambda^{n-1}+...+a_{n-1}\lambda+a_n=0,$$
then the function
$$f_k(X)=\exp(\lambda_kz)H_k(\uy)-\frac{1}{2\lambda_k}\exp(\lambda_k\overline{z})\puy H_k(\uy),$$
where $H_k$ is an arbitrary harmonic function of $\uy$, represents a solution of the differential equation $L_nf=0$. We can obtain a system of linearly independent solutions with different $\lambda_k$ by using an appropriate complete orthonormal system of harmonic functions $H_k(\uy)$.
\end{proposition}
Note that unlike homogeneous linear ordinary differential equations, where the maximum number of linearly independent solutions is equal to the order of the equation, in homogeneous linear partial differential equations this property is not satisfied. The set of solutions can have infinite dimension and not be limited by the order of the equation, as in the classical cases of Laplace, heat, or wave equations (see \cite{evans,gilbarg}). This is why the steering-type solutions we have found for \eqref{eqdh} are particular.

Let $\lambda_1=\lambda_2=...=\lambda_{\alpha_i}=0$ be a root of \eqref{eccara} of multiplicity $\alpha_i$ and all other roots be real simple roots.  Therefore, the first member of \eqref{eccara} has, in this case, the common factor $\lambda^{\alpha_i}$, that is, the coefficients $a_n=a_{n-1}=...=a_{n-\alpha_i+1}=0$. The corresponding homogeneous linear differential equation reduces to the following
\begin{equation}\label{eqmul}
a_0D^n f+a_1D^{n-1}f+...+a_{n-\alpha_i}D^{\alpha_i}f=0.
\end{equation}
A sequence of left monogenic polynomials $P_k(X)$ in a domain $\Omega\subset\R^{m+1}$ satisfying the recursion relation
\begin{itemize}
\item[(i)] $P_0(X)=1$,
\item[(ii)] $P_k(0)=0$,
\item[(iii)] $DP_k=kP_{k-1}$, $k=1,2,...$
\end{itemize}
is called a generalized Appell sequence (see \cite{GHS2}). We refer the reader to the preliminary work by Paul \'Emile Appell \cite{appell}. The property (iii) obviously generalizes the property of the “ordinary” powers $x^n$ or $z^n$. M.I. Falc\~ao, J.F. Cruz and H.R. Malonek were the first to construct the sequence
\begin{align*}
P_k(X)=\sum_{s=0}^kT_s^kX^{k-s}\overline{X}^s,\\
T_s^k=\binom{k}{s}\frac{\left(\frac{m+1}{2}\right)_{(k-s)}\left(\frac{m-1}{2}\right)_{(s)}}{m_{(k)}},
\end{align*}
where $a_{(k)}$ denotes the Pochhammer symbol, i.e. $a_{(k)}:=a(a+1)...(a+k-1)=\frac{\Gamma(a+k)}{\Gamma(a)}$, for any integer $k>1$, and $a_{(0)}:=1$ \cite{falcao1,malonek,falcao2}. It is clear that the sequence $\{P_k\}_{k=0}^{\alpha_i-1}$ solves equation \eqref{eqmul}. The property (i) can be altered, and in general we can consider any type of function, not just polynomials, so for example the sequence $\{z^kM(\uy)\}_{k=0}^{\alpha_{i}-1}$, where $\puy M(\uy)=0$, is also a generalized Appell sequence. This last sequence also solves equation \eqref{eqmul}. 

If the characteristic equation \eqref{eccara} has a root $\lambda_i\not=0$ of multiplicity $\alpha_i>1$, then the functions
\begin{equation}\label{71eq}
z^k\left[\exp(\lambda_iz)M(\uy)\right],\quad k=0,...,\alpha_i-1,
\end{equation}
with  $\puy M(\uy)=0$, are solutions to equation \eqref{eqdh}. The treatment of the case when equation \eqref{eccara} has a complex root is similar, with the important guideline that it is now appropriate to work on $\C_m:=\C\otimes\R_{0,m}$, where the imaginary unit $i\in\C$ commutes with the basic elements $e_j$, $j=1,...,m$.
 
Several generalizations of the exponential function have been discussed in Clifford analysis. This newly-defined exponential function $F(X)$ has a different form than the other generalizations reported in the literature. The exponential steering function $F(X)$ discussed here, in addition to satisfying the important relation \eqref{imDff}, is left monogenic. It is worth mentioning that in \cite[Definition 11.18, pag. 229]{GHS} the authors define an exponential function for paravectors of Clifford algebras, which is radially holomorphic. They introduce a Fueter transform to construct monogenic functions from their exponential function. The function \eqref{rel54} has the advantage of being monogenic, which is not the case for function $e^X$ defined in \cite{GHS}.

\section{Further results}
The Lam\'e–Navier system acting on vector fields of $\R^3$ is defined as follows
\begin{equation}
\mu\Delta\vec{u}+(\mu+\lambda)\nabla(\nabla\cdot\vec{u})=0,
\end{equation}
where $\mu$ and $\lambda$ are elastic constants commonly referred to as Lam\'e parameters. This system was originally introduced by the French mathematician Gabriel Lam\'e in 1837 while studying the method of separation of variables for solving the wave equation in elliptic coordinates \cite{lame}. A natural generalization of this system in a multidimensional setting can be conceived as
\begin{equation}\label{68}
\left(\frac{\mu+\lambda}{2}\right)\partial_X f\partial_X+\left(\frac{3\mu+\lambda}{2}\right)\partial_X^2f=0.
\end{equation}
We refer the reader to works \cite{DAS4,MD,MAB3} where this equation \eqref{68} is studied. 
Inframonogenic functions emerge as the solutions of the sandwich-type equation
\begin{equation}
\partial_X f\partial_X=0.
\end{equation}
Such functions were defined in \cite{MPS2}, where a Fischer decomposition for homogeneous polynomials in terms of inframonogenic polynomials was proved. Other results obtained for this functional class can be found in \cite{AbreuDAS,DAS1,DAS2,DAS3,DAS5,joao,Din,lavicka,MPS1,MAB1,MAB4,MMA,MDAS,wang}. Consider the exponential steering function
$$F(X)=\exp(z)A(\uy)+\exp(\overline{z})B(\uy).$$
It is easy to see that
\begin{align*}
\partial_XF(X)\partial_X&
=\exp(z)[\puy B+e_1\puy Be_1+\puy B\puy]+\exp(\overline{z})[\puy A-e_1\puy Ae_1+\puy A\puy++2B-2e_1Be_1+2B\puy].
\end{align*}
Therefore, if $\partial_X F(X)\partial_X=0$ then 
\begin{equation}\label{infras}
\left\{\begin{array}{rl}
\puy B+e_1\puy Be_1+\puy B\puy=0,\\
\puy A-e_1\puy Ae_1+\puy A\puy+2B-2e_1Be_1+2B\puy=0.
\end{array}\right.
\end{equation}
Let $I(\uy)$ be an inframonogenic function and $M(\uy)$ a right monogenic function, both on the variable $\uy$. The function 
\begin{equation}
\exp(z)[I(\uy)-e_1I(\uy)e_1]+\exp(\overline{z})[M(\uy)+e_1M(\uy)e_1]
\end{equation}
is inframonogenic. For example. if $I(\uy)=\frac{1}{2}(x_2^2+x_3^2)e_2e_4$ and $M(\uy)=\frac{1}{2}(x_2e_2-x_3e_3)$ then
$$\exp(z)(x_2^2+x_3^2)e_2e_4+\exp(\overline{z})(x_2e_2-x_3e_3)$$
is inframonogenic.

Using the calculations in section \ref{Section3}, the function $F(X)$ satisfies the generalized Lam\'e-Navier system \eqref{68} if
\begin{equation}\label{infras}
\left\{\begin{array}{rl}
\left(\frac{7\mu+3\lambda}{2}\right)\puy B+\left(\frac{\mu+\lambda}{2}\right)e_1\puy Be_1+\left(\frac{\mu+\lambda}{2}\right)\puy B\puy+\left(\frac{3\mu+\lambda}{2}\right)\puy^2 A=0,\\
\left(\frac{7\mu+3\lambda}{2}\right)\puy A-\left(\frac{\mu+\lambda}{2}\right)e_1\puy Ae_1+\left(\frac{\mu+\lambda}{2}\right)\puy A\puy+(7\mu+3\lambda)B-(\mu+\lambda)e_1Be_1+(\mu+\lambda)B\puy+\left(\frac{3\mu+\lambda}{2}\right)\puy^2B=0.
\end{array}\right.
\end{equation}
For instance, if $I(\uy)$ is as before, an inframonogenic function that depend only on the variable $\uy$, then 
$$\exp(z)[I(\uy)-e_1I(\uy)e_1]-\frac{1}{2}\exp(\overline{z})[\puy I(\uy)+e_1\puy I(\uy)e_1]$$
is solution to system \eqref{68} and does not depend on Lam\'e parameters. These solutions are commonly referred to as universal.

We can also construct functions of this type for so-called $(\alpha,\beta)$-monogenic functions, i.e. the solutions of system
\begin{equation}
\alpha f\partial_X+\beta\partial_X f=0,\quad\alpha,\beta\in\R.
\end{equation}
These functions were studied in \cite{RM}, where the authors found integral representation formulas related to the Lam\'e-Navier system. By direct calculation, we can see that $F(X)$ would be $(\alpha,\beta)$-monogenic if
\begin{equation}\label{eq7}
\left\{\begin{array}{rl}
\beta\puy B+\alpha A+\alpha e_1Ae_1+\alpha A\puy&=0,\\
\beta\puy A+(2\beta+\alpha)B-\alpha e_1Be_1+\alpha B\puy&=0.
\end{array}\right.
\end{equation}

It is also important to note that we can construct steering-type solutions for a wide  class of differential systems involving the generalized Cauchy-Riemann operator. For example,
$$\partial_X^2 F(X)\partial_X=0$$
$$\Updownarrow$$
$$\left\{\begin{array}{rl}
\puy^2 A+2\puy B+e_1\puy^2 Ae_1+2e_1\puy Be_1+\puy^2A\puy+2\puy B\puy=0,\\
2\puy A+\puy^2 B+4B-2e_1\puy Ae_1-e_1\puy^2Be_1-4e_1Be_1+2\puy A\puy+\puy^2B\puy+4B\puy=0,
\end{array}\right.$$
$$\partial_X F(X)\partial_X^2=0$$
$$\Updownarrow$$
$$\left\{\begin{array}{rl}
2\puy B+2e_1\puy Be_1+2\puy B\puy+\puy B\puy^2=0,\\
2\puy A-2e_1\puy Ae_1+2\puy A\puy +\puy A\puy^2+4B-4e_1Be_1+4B\puy+2B\puy^2=0.
\end{array}\right.$$
Nowadays, these functions are known as infrapolymonogenic and have been investigated in \cite{Minfra}. As a matter of fact, much interesting examples may be described.

\section*{Acknowledgments}
The first author wishes to express his gratitude to the financial support of the Postgraduate Study Fellowship of the  Secretar\'ia de
Ciencia, Humanidades, Tecnolog\'ia e Innovaci\'on (SECIHTI) (grant number 1043969). 

\section*{Declaration of Competing Interest}
The authors have no conflicts of interest to declare. 

\section*{Data availability}
Not applicable.

\section*{ORCID}
\noindent Daniel Alfonso Santiesteban: \url{https://orcid.org/0000-0003-0248-3942}\\
Dixan Pe\~na Pe\~na: \url{https://orcid.org/0000-0002-6588-5474}\\
Ricardo Abreu Blaya: \url{https://orcid.org/0000-0003-1453-7223}

\end{document}